\documentclass[preprint,12pt]{elsarticle}
\usepackage[utf8]{inputenc}
\usepackage{lineno,hyperref}
\modulolinenumbers[5]
\usepackage{graphicx}
\usepackage{pst-all}
\usepackage{pstricks-add}
\usepackage[crop=off]{auto-pst-pdf}

\usepackage{pst-pdf}
\usepackage[T1]{fontenc}
\usepackage{amsmath,amssymb,latexsym}
\usepackage{amsthm}
\usepackage[center,small,sf]{caption}
\usepackage{mathtools}
\usepackage{diagbox}

\newtheorem{defin}{Definition}

\newtheorem{prop}{Proposition}

\newcommand{\fg}{\mathfrak{g}}%

\usepackage[english]{babel}

\makeatletter
\def\ps@pprintTitle{%
   \let\@oddhead\@empty
   \let\@evenhead\@empty
   \let\@oddfoot\@empty
   \let\@evenfoot\@oddfoot
}
\makeatother
\begin{document}

\begin{frontmatter}

\title{
\textbf{Lie Algebra Classification, Conservation Laws and Invariant Solutions for a Variant of the Levinson-Smith Equation}\\
}

\author{Y. Acevedo \corref{cor1}}

\author{Danilo Hernández-García  
}
\author{Gabriel Loaiza}  
\author{Oscar Londoño}

\begin{abstract}
We obtain the optimal system's generating operators associated with the kind generalization of the Levinson–Smith equation. Using those operators we characterize all invariant solutions associated with this equation. 
Moreover, we present the variational symmetries and the corresponding conservation laws, using Noether's theorem. Finally, we classify the Lie algebra associated with the given equation.
\\
\end{abstract}


\end{frontmatter}

 \pagestyle{myheadings}
 \pagenumbering{arabic}

\section{Introduction}
Consider the following differential equation
\begin{equation}\label{ppp}
y_{xx}+\phi(x,y,y_x) y_x=\gamma(x,y),
\end{equation}
where $\gamma$ is non-linear function, whose meaning "function of disturbance" and the function $\phi$ can be identified as a non-linear coefficient of friction. This type of equations \eqref{ppp} are known as variants of the Levinson–Smith equation. In \cite{Kamke1977}, Kamke considers a particular case of \eqref{ppp}
\begin{equation}\label{nn1}
y_{xx}=x^{-a}\left(xy_x-y\right)^{a} f(x,y)=0,
\end{equation}
where $a$ is real number and $f$ is arbitrary function, in this work,  Kamke proposes a solution in the case $a=2$ and $f(x,y)=-(x+y)^{-1}$, this solucion is 
\begin{equation}\label{nn2}
y(x)=-x+x C_{1} e^{ \frac{C_{2}}{x}},\,\,\mbox{where $C_2,C_1$ are constants.}
\end{equation}
In \cite{Polyanin2002}, Polyanin and Zaitsev introduce a particular case of \eqref{ppp}
\begin{equation}\label{nn3}
(y+ax) y_{x x}=b x^{n}\left(x y_{x}-y\right)^{2}
\end{equation}
where $a,b,n$ are real numbers, and in this work a substitution is proposed
$y=-a x+x z$ leads to the equation $x z z_{r r}+2 z z_{r}-b x^{n+3}\left(z_{r}\right)^{2}=0.$ Having set $w=\frac{z_{x}}{z},$ we obtain a Bernoulli equation $x w_{r}+2 w+x\left(1-b x^{n+2}\right) w^{2}=0.$ In \cite{Duarte1998}, Duarte et al introduce a particular case of \eqref{ppp}, which in turn, is a particularization of \eqref{nn1} and \eqref{nn3}
 \begin{equation}\label{ODE}
y_{xx}=-y_x^2(x+y)^{-1}-2yy_xx^{-1}(x+y)^{-1}-y^2x^{-2}(x+y)^{-1},
\end{equation}
and its solution, 
\begin{equation}\label{qw}
y(x)=x e^{\left(-\frac{C_{1}}{x}+C_{2}\right)}-x,\,\,\mbox{where $C_2,C_1$ are constants.}
\end{equation}
In this work, they present the Lie group of symmetries of \eqref{ODE}, using a ODEtools Maple package. \\
In \cite{Hindawi2021}, another type of Levinson-Smith equation is proposed, which is:
$$y_{x x}=2 y_{x}^{2} y^{-1}+x^{-1} y_{x}+x^{-1} y^{2},$$
for this equation, the authors presented various implicit and explicit solutions, these solutions are calculated using optimal algebra.\\

The proposal of our work is: $i)$ to calculate the $3-$dimensional Lie symmetry group in all detail, $ii)$ to present the optimal system (optimal algebra) for \eqref{ODE}, $iii)$ making use of all elements of the optimal algebra, to propose invariant solutions for \eqref{ODE}, then  $iv)$ to construct the Lagrangian with which we could determine the variational symmetries using Noether${}'$s theorem and build some non-trivial conservation laws, and  finally $v)$ to obtain a representation of the Lie algebra corresponding to the Lie group symmetries for the equation \eqref{ODE}.

\section{Lie Point symmetries}

First we explore the Lie symmetry group for \eqref{ODE}. The main result of this section is:
\begin{prop}\label{PROP}
The Lie point symmetries for the equation \eqref{ODE} are generated by the vector fields:
\begin{align}\label{Simmetrias}
    &\Pi_1=x\frac{\partial}{\partial x}+y\frac{\partial}{\partial y},& &\Pi_2=-x\frac{\partial}{\partial x}+x\frac{\partial}{\partial y},& &\mbox{and}& &\Pi_3=x^2\frac{\partial}{\partial x}+xy\frac{\partial}{\partial y}.&\\ \nonumber
\end{align}
\end{prop}

\begin{proof}[Proof] 
The general form for the generator operators of a Lie group of a parameter admitted by \eqref{ODE} is:
$$x\rightarrow x+\epsilon \xi(x,y)+O\left(\epsilon^{2}\right),\,\,\,\,\,\,\,\,\,\,\,\,\text{and}\,\,\,\,\,\,\,\,\,\,\,\,\,\,y\rightarrow y+\epsilon \eta(x,y)+O\left(\epsilon^{2}\right),$$
where $\epsilon$ is the group parameter. The vector field associated with this group of transformations is
$
\Gamma=\xi(x,y)\frac{\partial}{\partial x}+\eta(x,y)\frac{\partial}{\partial y},
$ with $\xi,\eta$ 
differentiable functions in $\mathbb{R}^2$. 
To find the infinitesimals $\xi(x,y)$ and $\eta(x,y)$, 
we applied the second extension operator,
\begin{eqnarray}\label{eq4}
\Gamma^{(2)} &=&\Gamma +\eta_{[x]}\frac{\partial}{\partial {y_x}}+\eta_{[xx]}\frac{\partial}{\partial {y_{xx}}},
\end{eqnarray}
to the equation  \eqref{ODE}, 
obtaining the following symmetry condition
\begin{eqnarray}\label{3eq6}
 \xi(x+y)^{-1}\left[-y_{x}^2(x+y)^{-1}-2yy_{x}x^{-1}(x+y)^{-1}-2yy_{x}x^{-2}-y^2x^{-2}(x+y)^{-1}\right.&\nonumber \\
 \left.-2y^2x^{-3}\right]+\eta(x+y)^{-1}\left[-y_{x}^2(x+y)^{-1}-2yy_{x}x^{-1}(x+y)^{-1}+2y_{x}x^{-1} \right.&\nonumber \\
 \left.-y^2x^{-2}(x+y)^{-1}+2yx^{-2} \right]+\left(2y_{x}-2yx^{-1}\right)(x+y)^{-1}\eta_{[x]}+\eta_{[xx]}=0.&
\end{eqnarray}
 
where  $\eta_{[x]}$ and $\eta_{[xx]}$ 
are the coefficients in $\Gamma^{(2)}$ 
given by:
\begin{eqnarray}\label{3eq4.1}
\eta_{[x]}&=&D_x[\eta]-(D_x[\xi])y_x=\eta_x +(\eta_y-\xi_x)y_x-\xi_yy_x^2. \nonumber\\
\eta_{[xx]}&=&D_x[\eta_{[x]}]-(D_x[\xi])y_{xx}, \nonumber \\
           &=&\eta_{xx}+(2\eta_{xy}-\xi_{xx})y_x+(\eta_{yy}-2\xi_{xy})y_x^2-\xi_{yy}y_x^3\nonumber\\
           & &+(\eta_y-2\xi_x)y_{xx}-3\xi_yy_xy_{xx}.
\end{eqnarray}
where $D_{x}$ 
is the total derivative operator:
$
D_{x}=\partial_{x}+y_{x}\partial_{y}+y_{xx}\partial_{y_{x}}+\cdots  $.

After applying \eqref{3eq4.1} in \eqref{3eq6} and 
substitute in the resulting expression $y_{xx}$ by \eqref{ODE}, 
is obtained:
\begin{eqnarray}\label{D1}
[\xi_y(x+y)^{-1}-\xi_{yy}]y_{x}^3+[-\xi(x+y)^{-2}-\eta(x+y)^{-2}+\eta_{y}(x+y)^{-1}+\eta_{yy}&\nonumber \\
-2\xi_{xy}+8yx^{-1}\xi_{y}(x+y)^{-1}]y_{x}^2+[-2yx^{-1}\xi(x+y)^{-2}-2yx^{-2}\xi(x+y)^{-1}&\nonumber \\
-2yx^{-1}\eta(x+y)^{-2}+2x^{-1}\eta(x+y)^{-1}+2\eta_x(x+y)^{-1}-4yx^{-1}\eta_{y}(x+y)^{-1}& \nonumber \\
+2yx^{-1}\xi_x(x+y)^{-1}+4y\xi_{x}(x+y)^{-1}+3y^2x^{-2}\xi_{y}(x+y)^{-1}+2\eta_{xy}-\xi_{xx}]y_x+& \nonumber \\
\left[-y^2x^{-2}\xi(x+y)^{-2}-2y^2x^{-3}\xi(x+y)^{-1}-y^2x^{-2}\eta(x+y)^{-2}+2yx^{-2}\eta(x+y)^{-1}\right.& \nonumber \\
\left. -2yx^{-1}\eta_{x}(x+y)^{-1}-y^2x^{-2}\eta_{y}(x+y)^{-1}+2y^2x^{-2}\xi_{x}(x+y)^{-1}+\eta_{xx}\right]=0.&\nonumber
\end{eqnarray}
Then, analyzing the coefficients with respect to the independent variables $y_{x}^3,y_{x}^2,y_{x},1$ we get the following system of determining equations:
\begin{subequations}
\begin{align}
 \xi_y-\xi_{yy}(x+y)=0, \label{Det1} \\
  -\xi-\eta+\eta_{y}(x+y)+(\eta_{yy}-2\xi_{xy})(x+y)^2+8yx^{-1}\xi_{y}(x+y)=0,& \label{Det2}\\
  -2yx^{-1}\xi-2yx^{-2}\xi(x+y)-2yx^{-1}\eta+2x^{-1}\eta(x+y)+2\eta_{x}(x+y)&\nonumber\\
  -4yx^{-1}\eta_{y}(x+y)+2yx^{-1}\xi_{x}(x+y)+4y\xi_{x}(x+y)+3y^2x^{-2}\xi_{y}(x+y) &\nonumber\\
  +(2\eta_{xy}-\xi_{xx})(x+y)^2=0,& \label{Det3}\\
    -y^2x^{-2}\xi-2y^2x^{-3}\xi(x+y)-y^2x^{-2}\eta+2yx^{-2}\eta(x+y)-2yx^{-1}\eta_{x}\nonumber\\
   (x+y)-y^2x^{-2}\eta_{y}(x+y)+2y^2x^{-2}\xi_{x}(x+y)+\eta_{xx}(x+y)^2=0.&\label{Det4} 
\end{align}
\end{subequations}
Solving \eqref{Det1} - \eqref{Det4} we get $\xi=k_1x-k_2x+k_{3}x^2$ and $\eta=k_{1}y+k_{2}x+k_{3}xy$ where $k_1,k_2,k_3,k_4$ and $k_5$ are arbitrary constants. Therefore, the generators of the Lie point symmetries of \eqref{ODE} are the operators $\Pi_1$ - $\Pi_3$ outlined in the declaration of the Proposition \ref{PROP}.
\end{proof}

\section{Optimal Algebra}
In this section, we take into account the methodology presented in \cite{Olver1986,hydon2000symmetry,Zahid2009,Zewdie}, we obtain optimal algebra associated to the symmetry group of \eqref{ODE}. 

First, we should calculate the commutator table, which is obtained from the operator
\begin{equation}\label{n1}
[\Pi_{\alpha},\Pi_{\beta}]=\Pi_\alpha \Pi_\beta-\Pi_\beta \Pi_\alpha=\sum_{i=1}^n\left(\Pi_\alpha(\xi_{\beta}^i)-\Pi_\beta(\xi_{\alpha}^i) \right) \frac{\partial}{\partial x^i},
\end{equation}
where $i=1,2$, with $\alpha, \beta=1,\cdots,3$  and  $\xi_{\alpha}^i,\xi_{\beta}^i$ are the respective coefficients of the  operators $\Pi_{\alpha},\Pi_{\beta}.$ After applying the operator \eqref{n1} to the symmetry group of \eqref{ODE}, we obtain the Table \ref{tabla ODE}.
\begin{table}
\centering
\resizebox{4cm}{!} {
\begin{tabular}{||c||c|c|c||}
\hline
\hline
& $\Pi_1$ & $\Pi_2$ & $\Pi_3$   \\ \hline \hline
$\Pi_1$ & $0$ & $0$ & $\Pi_3$  \\ \hline
$\Pi_2$ & $0$ & $0$ & $-\Pi_3$  \\ \hline
$\Pi_3$ & $-\Pi_3$ & $\Pi_3$ & $0$  \\  \hline \hline 
\end{tabular}
}
\caption{Commutators table associated to the symmetry group of \eqref{ODE}.}\label{tabla ODE}
\end{table}
Now, we calculate the adjoint action representation of the symmetries of \eqref{ODE} using the Table \ref{tabla ODE} and the operator 
\begin{equation*}\label{Op.autoadjunto}
Ad(exp( \lambda \Pi))G'=\sum_{n=0}^{\infty} \frac{\lambda^n}{n!}(ad(\Pi))^nG\,\,\,\text{for the symmetries} \,\,G\, \text{and}\,\Pi.
\end{equation*}
The Table \ref{tabla Adj} shows the adjoint representation obtained for each $\Pi_i$ using the previous operator.
\begin{table}
\centering
\resizebox{6cm}{!} {
\begin{tabular}{||c||c|c|c||}
\hline
\hline
Adj[\,,\,]& $\Pi_1$ & $\Pi_2$ & $\Pi_3$ \\ \hline \hline
$\Pi_1$ & $\Pi_1$ & $\Pi_2$ & $e^{-\lambda}\Pi_3$     \\ \hline
$\Pi_2$ & $\Pi_1$ &$\Pi_2$ & $e^{\lambda}\Pi_3$   \\ \hline
$\Pi_3$& $\Pi_1+\lambda \Pi_3$ &$\Pi_2-\lambda \Pi_3$ & $\Pi_3$  \\  \hline \hline
\end{tabular}
}
\caption{Adjoint representation table for the Lie symmetries of \eqref{ODE}.}\label{tabla Adj}
\end{table}

\begin{prop}\label{uno}
The equation \eqref{ODE} has the optimal algebra is given by
$$\Pi_1+\Pi_2, \Pi_1+b_2\Pi_3,$$
$$a_1\Pi_1+a_2\Pi_2 \,,\,\text{with}\,\,a_1\neq a_2,$$
$$a_1\Pi_1+a_1\Pi_2+\Pi_3,$$
$$a_1\Pi_1+a_2\Pi_2+\Pi_3\,,\,\text{with}\,\,a_1\neq 1,$$
\end{prop}

\begin{proof}
To calculate the optimal algebra system, we start with the generators of symmetries \eqref{Simmetrias} and a generic nonzero vector. Let 
\begin{equation}\label{y1}
G=a_{1} \Pi_{1}+a_{2} \Pi_{2}+a_{3} \Pi_{3}.
\end{equation}
The objective is to simplify as many coefficients $a_i$ as possible, through maps adjoint  to $G$, using Table \ref{tabla Adj}.
\begin{enumerate}
   \item[1)] Assuming $a_3 = 1$ in \eqref{y1} we have that $G = a_{1} \Pi_{1}+a_{2} \Pi_{2}+\Pi_{3}$. Applying the adjoint operator to $(\Pi_1,G)$ and $(\Pi_2,G)$ we don't have any reduction, on the other hand applying the adjoint operator to $(\Pi_3,G)$ we get
\begin{equation}\label{y2}
G_{1}=A d\left(\exp \left(\lambda_{1} \Pi_{3}\right)\right) G=a_1\Pi_1+a_2\Pi_2+(1+\lambda_{1}(a_1-a_2))\Pi_3.
\end{equation}

\textbf{$1.1)$ Case $a_1 -a_2\neq 0$.} Using $\lambda_1=\frac{1}{a_1-a_2}$ with $a_1\neq a_2$, in \eqref{y2}, $\Pi_3$ is eliminated, therefore $G_1 =a_1\Pi_1+a_2\Pi_2$. Thus, we obtain the first element
\begin{equation}\label{optimal1}
G_1=a_1\Pi_1+a_2\Pi_2,\,\,\mbox{with $a_1\neq a_2$.}  
\end{equation}
 So, the first reduction of the generic element \eqref{y1} was done.\\
\textbf{$1.2)$ Case $a_1=a_2$.} Thus $G_1=a_1\Pi_1+a_1\Pi_2+\Pi_3$. We can not reduce more, then we have other element ot the optimal algebra
\begin{equation}\label{optimal2}
G_1=a_1\Pi_1+a_1\Pi_2+\Pi_3.
\end{equation}
So, the first reduction of the generic element \eqref{y1} was done.\\
\item[2)] Assuming $a_3 = 0$ and $a_2=1$ in \eqref{y1}, we have that $G = a_{1} \Pi_{1}+ \Pi_{2}$. Applying the adjoint operator to $(\Pi_1,G)$ and $(\Pi_2,G)$ we don't have any reduction, on the other hand applying the adjoint ope\-rator to $(\Pi_3,G)$ we get
\begin{equation}\label{y22W}
G_{2}=A d\left(\exp \left(\lambda_{2} \Pi_{3}\right)\right) G= a_1\Pi_1+\Pi_{2}+\lambda_2(a_1-1)\Pi_{3}.
\end{equation}

\textbf{$2.1)$ Case $a_1-1 \neq 0$.} Using $\lambda_2=\frac{b_1}{a_1-1}$, with $a_1\neq 1$, in \eqref{y22W}, we have other element of the optimal algebra 
\begin{equation}\label{elemt1K}
G_2=a_1\Pi_1+\Pi_2+b_1\Pi_3,\,\,\mbox{with $a_1\neq1$.}
\end{equation}
Then, other reduction of the generic element \eqref{y1} ends.\\
\textbf{$2.2)$ Case $a_1=1$.} Then $G_2 =\Pi_{1}+\Pi_2.$ Thus, we get other member of the optimal algebra 
\begin{equation}\label{elemt1KZZ}
G_2=\Pi_1+\Pi_2.
\end{equation}
So, other reduction of the generic component \eqref{y1} ends.\\
\item[3)] Assuming $a_3= a_2 = 0$ and $a_1=1$ in \eqref{y1}, we have that $G = \Pi_{1}$. Applying the adjoint operator to $(\Pi_1,G)$ and $(\Pi_2,G)$ we don't have any reduction, on the other hand applying the adjoint ope\-rator to $(\Pi_3,G)$ we get
\begin{equation}\label{y22}
G_{3}=A d\left(\exp \left(\lambda_{3} \Pi_{3}\right)\right) G= \Pi_1+\lambda_{3}\Pi_3.
\end{equation}
We don't get any reduction, thus using $\lambda_{3}=b_2$ we have other element of the optimal algebra
\begin{equation}\label{elemtKKL1}
G_3=\Pi_1+b_2\Pi_3.  
\end{equation}
So, other reduction of the generic component \eqref{y1} ends.

\end{enumerate}
\end{proof}

\section{Invariant solutions by  generators of the optimal algebra}
In this section, we characterize all invariant solutions taking into account some operators that generate the optimal algebra presented in Proposition \ref{uno}. For this purpose, we use the method of invariant curve condition \cite{hydon2000symmetry} (presented in section 4.3), which is given by the following equation
 \begin{equation}\label{con1}   \mathit{Q}(x,y,y_x)=\eta-y_x\xi=0.
\end{equation}
Using the element $\Pi_1+\Pi_2$ from Proposition \ref{uno}, under the condition \eqref{con1}, we obtain that $\mathit{Q}=\eta_{1+2}-y_x\xi_{1+2}=0,$  which implies $(y+x)-y_x(0)=0$, then solving this ODE, we have $y(x)=-x$, which is a solution explicit for \eqref{ODE}, but it makes a null denominator too. Using an analogous procedure with all of the elements of the optimal algebra (Proposition \ref{uno}, we obtain both implicit and explicit invariant solutions that are shown in the Table \ref{CompSolutions}, with $c$ being a constant.
\begin{table}
\centering
\resizebox{\textwidth}{!}{%
\begin{tabular}{|c|c|c|c|c|}
\hline
&\textbf{Elements} & \textbf{\textbf{$Q(x,y,y_x)=0$}} & \textbf{Solutions} & \textbf{Type Solution} \\ \hline \hline
$1$ &$\Pi_1+\Pi_2$ & $(y+x)-y_x(0)=0$ & $y(x)=-x$ & 
Explicit, but it makes a null denominator too. \\ \hline
$2$&$\Pi_1+\Pi_3$ & $(xy+y)-y_x(x^2+x)=0$ & $y(x)=cx$ & Explicit \\ \hline
$3$&$-\Pi_1+\Pi_2$ & $(x-y)-y_x(-2x)=0$ & $y(x)=c \sqrt{x}-x$ & Explicit \\ \hline
$4$&$\Pi_1+\Pi_2+\Pi_3$ & $(xy+y+x)-y_x(x^2)=0$ & $y(x)=c e^{-1 / x} x-x$  & Explicit \\ \hline
$5$&$-\Pi_1+\Pi_2+\Pi_3 $ & $(xy-y+x)-y_x(x^2-2x)=0$ & $y(x)=c \sqrt{2-x} \sqrt{x}-\frac{\sqrt{x} \sqrt{-(x-2) x}}{\sqrt{2-x}}$ & Explicit  \\ \hline
\hline
\end{tabular}%
}
\caption{Solutions for \eqref{ODE} using invariant curve condition.} \label{CompSolutions}
\end{table}
$\mathbf{Remarks:}$ Note that the solution in the numeral $4$ is a particular case of the solution \eqref{qw} and \eqref{nn2} presented in \cite{Duarte1998} and \cite{Kamke1948} .
\section{Variational symmetries and conserved quantities}
In this section, we present the variational symmetries of \eqref{ODE} and we are going to use them to define conservation laws via Noether's theorem \cite{Noether1918}. First of all, we are going to determine the Lagrangian using the Jacobi Last Multiplier method, presented by Nucci in \cite{Nucci2009}, and for this reason, we are urged to calculate the inverse of the determinant $\Delta$,   
$$
\Delta = \begin{vmatrix}
 x&  y_x& y_{xx}\\ 
 \Pi_{1,x}&  \Pi_{1,y}& \Pi_{1}^{(1)}\\ 
 \Pi_{2,x}&  \Pi_{2,y}& \Pi_{2}^{(1)}
\end{vmatrix}= \begin{vmatrix}
 x&  y_x& y_{xx}\\ 
 x&  y& 0\\ 
 -x&  x& 1-y_x
\end{vmatrix},
$$
where $\Pi_{1,x},\Pi_{1,y},\Pi_{2,x},$ and $\Pi_{2,y}$ are the components of the symmetries $\Pi_1,\Pi_2$ shown in the Proposition \ref{Simmetrias}  and $\Pi_{1}^{(1)},\Pi_{2}^{(1)}$  as its first prolongations. Then we get $\Delta =x(x+2y-y_x)$ which implies that $M=\frac{1}{\Delta}=\frac{x^{-1}}{x+2y-y_x}.$ Now, from \cite{Nucci2009}, we know that $M$ can also be written as $M=L_{y_xy_x}$ which means that $L_{y_xy_x}=\frac{x^{-1}}{x+2y-y_x},$ then integrating twice with respect to $y_x$ we obtain the Lagrangian
\begin{eqnarray}\label{lan}
L(x,y,y_x)&=& -x^{-1}y_x\ln (x+2y-y_x)+(1+2x^{-1}y)\ln (x+2y-y_x)\nonumber\\
&+&x^{-1}y_x+y_xf_1(x,y) + f_2(x,y),
\end{eqnarray}
where $f_1,f_2$ are arbitrary functions. From the preceding expression we can consider $f_1=f_2=0.$ It's possible to find more Lagrangians for
\eqref{ODE} by considering other vector fields given in the Proposition \ref{Simmetrias}. We then calculate
\begin{equation}\label{oo1}
    \xi(x,y) L_x+\xi_x(x,y) L+\eta(x,y) L_y +\eta_{[x]}(x,y)L_{y_x}=D_x[f(x,y)],
\end{equation}
using \eqref{lan} and \eqref{3eq4.1} into \eqref{oo1} and rearranging and associating terms with respect to $1,y_x,y_x^2,\frac{1}{x+2y-y_x},\frac{y_x}{x+2y-y_x},\frac{y_x^2}{x+2y-y_x},\frac{y_x^3}{x+2y-y_x},\ln (x+2y-y_x),y_x\ln (x+2y-y_x)$ and $y_x^2\ln (x+2y-y_x),$ we obtain the following determining equations, presented in the following table
\begin{table} 
\centering
\resizebox{\textwidth}{!}{%
\begin{tabular}{|c|c|c|}
\hline
$\xi_y=0$ & \textbf{$\eta_y-\xi_x+(x+2y)\xi_y=0$} & \textbf{$x^{-2}\xi-x^{-1}\eta_y=0$} \\ \hline
\textbf{$\xi+\eta-\eta_x=0$} & $-x^{-2}\xi+x^{-1}\eta_y-f_y=0$ & $x^{-1}\eta_x-f_x=0$ \\ \hline
\textbf{$-\xi-2\eta+\eta_x-(x+2y)(\eta_y-\xi_x)=0$} & $x^{-1}(\xi+2\eta-\eta_x)-(x^{-1}+2yx^{-2})\xi+(1+2yx^{-1})\xi_x=0$ &  \\ \hline
\end{tabular}%
}\caption{Determinant equations} \label{tab1}
\end{table}

Solving the preceding system (Table \ref{tab1} for $\xi,\eta$ and $f$ we obtain the infinitesimal generators of Noether's symmetries
\begin{align}\label{b}
    &\eta=a_2e^{2x},& &\xi=0& &\mbox{and}&  &f(x)=a_2xe^{2x}-\frac{a_2e^{2x}}{2}.&
\end{align}
where $a_2$ is arbitrary constant. Then, the Noether symmetry group or variational symmetry is
\begin{align}\label{b1}
    &V_1=e^{2x}\frac{\partial}{\partial y}.&
\end{align}
According to \cite{Fomin2000}, in order to obtain the conserved quantities or conservation laws, we should solve 
$$
I=(X y_x-Y)L_{y_x}-X L+f,
$$
so, using \eqref{lan}, \eqref{b} and \eqref{b1}. Therefore, the conserved quantity is given by
\begin{eqnarray}\label{b2}
I_1&=&\frac{e^{2x}x^{-1}y_x}{((x+2y)-y_x)}-e^{2x}x^{-1}\ln (x+2y-y_x)-\frac{e^{2x}(1+2x^{-1}y)}{((x+2y)-y_x)}\nonumber\\
&+&e^{2x}x^{-1}+a_2xe^{2x}-\frac{a_2e^{2x}}{2},
\end{eqnarray}
    \section{Classification of  Lie algebra}
Usually, a finite dimensional Lie algebra in a field of characteristic $0$ is classified by the Levi's theorem. It is, there exist  two important classes of Lie algebras, the solvable and the semisimple. In each classes, there are some particular classes that get other classification, for example, in the solvable one, we get the nilpotent Lie algebra.\\

According to the  Table \ref{tabla ODE}, we have a three dimensional Lie algebra, to corroborate it, first, we should be remember some basic criteria to classify a  Lie algebra for solvable and semisimple Lie algebra.
We will denote $K(.,.)$ the Cartan-Killing form. Some of the following propositions can be found in \cite{humphreys2012introduction}.
\begin{prop}\label{semisimple2}
{(Cartan's theorem)} A Lie algebra is semisimple if and only if its Killing form is nondegenerate.
\end{prop}
\begin{prop}\label{solbable1}
A Lie subalgebra $\mathfrak{g}$ is solvable if and only if $K(X,Y)=0$ for all $X\in [\mathfrak{g},\mathfrak{g}]$ and $Y\in \mathfrak{g}$. Other way to write that  is $K(\mathfrak{g},[\mathfrak{g},\mathfrak{g}])=0$.
\end{prop}
We also need  the next statements to make the  classification.
\begin{defin}
Let $\fg$ be a finite-dimensional Lie algebra over an arbitrary field $k$. Choose a basis ${e_{j}}$,
$1\leq i \leq n$, in $\fg$ where $n= \text{dim } \fg $ and set $[e_i,e_j]=C_{ij}^{k}e_{k}$, then the coefficients $C_{ij}^{k}$ are named structure constants.
\end{defin}
\begin{prop}\label{propo1}
Let $\fg_{1}$ and $\fg_{2}$ be two Lie algebras of dimension $n$. Suppose each has a basis with respect to which the structure constant are the same, then $\fg_{1}$ and $\fg_{2}$ are isomorphic.
\end{prop}

Let $\mathfrak{g}$ the Lie algebra related to the symmetry group of infinitesimal generators of the equation (1)
as stated by the table of the commutators, it is enough to consider
the next relations: $[\Pi_1,\Pi_3]=\Pi_3$. Now, we calculate Cartan-Killing form $K$ as follows.
$$K=\left[\begin{matrix}  
	 1 & 0 &  0  \\
     0 & 0 &  0  \\
     0 & 0 &  0  
	\end{matrix} \right],$$
	which the determinant vanishes, and thus by Cartan criterion it is not semisimple, (see Proposition \ref{semisimple2}. Since a nilpotent Lie algebra has a Cartan-Killing form that is identically zero, we  conclude, using the counter-reciprocal of the last claim, that the Lie algebra $\mathfrak{g}$ is not nilpotent.\\
	
	We verify that the Lie algebra is solvable using the Cartan criteria to solvability, (Proposition \ref{solbable1}, and then we have a   solvable nonnilpotent  Lie algebra.
	The Nilradical of the Lie algebra $\mathfrak{g}$ is generated by $\Pi_2, \Pi_3$,  and his center is generated by $\Pi_2$. Then, we have a solvable nonnilpotent Lie algebra. Furthermore, If we make $e_1:=-\Pi_1$ and $e_2:=-\Pi_3$ we get 
	$[e_1,e_2]=e_1$ from that,  and doing the comparisom with the Bianchi's classification of three dimensional Lie algebras and   by the Proposition \ref{propo1} we have a isomorphism of the  Lie algebra $\mathfrak{g}$  with  $\mathfrak{g}_{2,1}\oplus \mathfrak{g}_{1}$, where  $\mathfrak{g}_{2,1}=\mathfrak{af}(1)$,
	ie., according the Bianchi's classification we have a
	a decomposable Solvable, Lie algebra  Bianchi of   type III, see \cite{bianchi1897sugli}.
	
  In summery we have the next proposition.
  \begin{prop}
The  $3$-dimensional Lie algebra $\mathfrak{g}$ related to the symmetry group of the equation (1) is a solvable nonnilpotent. Moreover, the Lie algebra is isomorphich with $\mathfrak{g}_{2,1}\oplus \mathfrak{g}_{1}$  in the Bianchi's classification, (decomposable Solvable, Lie algebra  Bianchi of   type III) .
\end{prop}
\section{Conclusion}
Using the Lie symmetry group (see Proposition \eqref{PROP}, we calculated the optimal algebra, as it was presented in Proposition \eqref{uno}. Using these operators it was possible to characterize all the invariant solutions (see Table \eqref{CompSolutions}, the solution in the numeral $4$ is a particular case of the solution \eqref{qw} and \eqref{nn2} presented in \cite{Duarte1998} and \cite{Kamke1948}, the rest of these solutions do not appear in the literature known until today.

It has been shown the variational symmetry for \eqref{ODE} in \eqref{b1} with its corresponding conservation law \eqref{b2}. The Lie algebra associated to the equation \eqref{ODE} decomposable Solvable, Bianchi type III. Therefore, the goal initially proposed was achieved.

For future works, the equivalence group theory could be also considered to obtain preliminary classifications associated to a complete classification of \eqref{ODE}.
\section*{Acknowledgments}

Y. Acevedo, G. Loaiza and O.M.L Duque are grateful to EAFIT University, Colombia, for the financial support in the project "Study and applications of diffusion processes of importance in health and computation" with code $11740052022$.\\

\bibliography{references.bib}

\end{document}